\documentclass[12pt]{amsart}
\usepackage{amssymb,hyperref}

\newtheorem{theorem}{Theorem}[section]
\newtheorem{definition}[theorem]{Definition}

\newtheorem{corollary}[theorem]{Corollary}

\newtheorem{question}[theorem]{Question}

\begin{document}

\title{The orthogonal Weingarten formula in compact form}

\author{Teodor Banica}
\address{T.B.: Department of Mathematics, Cergy-Pontoise University, 2 avenue Chauvin, 95302 Cergy-Pontoise, France. {\tt teodor.banica@u-cergy.fr}}


\begin{abstract}
We present a compact formulation of the orthogonal Weingarten formula, with the traditional quantity $I(i_1,\ldots,i_{2k}:j_1,\ldots,j_{2k})=\int_{O_n}u_{i_1j_1}\ldots u_{i_{2k}j_{2k}}\,du$ replaced by the more advanced quantity $I(a)=\int_{O_n}\Pi u_{ij}^{a_{ij}}\,du$, depending on a matrix of exponents $a\in M_n(\mathbb N)$. Among consequences, we establish a number of basic facts regarding the integrals $I(a)$: vanishing conditions, possible poles, asymptotic behavior.
\end{abstract}

\maketitle

\section*{Introduction}

The computation of polynomial integrals over the orthogonal group $O_n$ is a key problem in mathematical physics. These integrals are indeed known to appear in a wealth of concrete situations, coming from random matrices, lattice models, combinatorics.

The standard approach to the computation of such integrals is via the Weingarten formula. This formula, originating from Weingarten's paper \cite{wei}, and worked out by Collins in \cite{col}, then by Collins and \'Sniady in \cite{csn}, is an identity of the following type:
$$\int_{O_n}u_{i_1j_1}\ldots u_{i_{2k}j_{2k}}\,du=\sum_{\pi,\sigma}\delta_\pi(i)\delta_\sigma(j)W_{kn}(\pi,\sigma)$$

Here the sum is over all pairings of $\{1,\ldots,2k\}$, also called Brauer diagrams \cite{bra}, and the delta symbols, describing the coupling between indices and diagrams, are 0 or 1. As for $W_{kn}$, this is the key combinatorial ingredient: the Weingarten function.

The exact or asymptotic computation of $W_{kn}$ is a quite subtle problem, and several results have been recently obtained on this subject. Let us mention here the work of Collins and Matsumoto \cite{cma}, Matsumoto and Novak \cite{mno} and Zinn-Justin \cite{zin}, providing a deep insight into the combinatorics of $W_{kn}$. Also, the foundational part of theory has benefited from several abstract versions and generalizations, developed in \cite{bco}, \cite{cur}, \cite{dif}.

The starting point for the considerations in the present paper is the following elementary reformulation of the Weingarten formula:
$$\int_{O_n}\prod_{i=1}^n\prod_{j=1}^nu_{ij}^{a_{ij}}\,du=\sum_{\pi,\sigma}\delta_\pi(a_l)\delta_\sigma(a_r)W_{kn}(\pi,\sigma)$$

Here $a\in M_n(\mathbb N)$ is the matrix of exponents appearing in the original Weingarten formula, and the delta symbols are once again 0 or 1. This formula is of course fully equivalent to the original Weingarten one, but the slight difference in the formulation leads to some potentially interesting consequences, that we will explore in this paper.

The idea is that our matrix formulation, with $a\in M_n(\mathbb N)$ replacing the multi-index $(i_1,\ldots,i_{2k}:j_1,\ldots,j_{2k})\in\mathbb\{1,\ldots,n\}^{4k}$, allows us to escape a bit from the heavy combinatorial context of \cite{wei}, \cite{col}, \cite{csn}, \cite{cma}, \cite{mno}, \cite{zin}. As an example, in the case where $a$ has only a nonzero row, with even entries, the exact formula for our integral is as follows:
$$\int_{O_n}u_{11}^{a_1}\ldots u_{1n}^{a_n}=\frac{(n-1)!!a_1!!\ldots a_n!!}{(n+\Sigma a_i-1)!!}$$

This kind of result, while very basic, is barely visible with the original Weingarten formula, and provides an excellent motivation for our proposed $a\in M_n(\mathbb N)$ reformulation. Some other exact formulae of this type are actually available, see the comments below.

So, let us go back now to the modified Weingarten formula, with the matrix $a\in M_n(\mathbb N)$ replacing the multi-indices $i,j$. This provides a combinatorial expression for the following integral, to be regarded from now on as being the ``basic quantity'' to be computed:
$$I(a)=\int_{O_n}\prod_{i=1}^n\prod_{j=1}^nu_{ij}^{a_{ij}}\,du$$

In this paper we will browse through the various combinatorial formulae in \cite{wei}, \cite{col}, \cite{csn}, \cite{cma}, \cite{mno}, \cite{cur}, by converting them into concrete results about $I(a)$. We will obtain in particular exact results regarding the vanishing, sign and poles of $I(a)$, viewed as a rational function of $n\in\mathbb N$, and also about its asymptotic behavior, with $n\to\infty$.

We believe that all these results can be of help in connection with a key problem, namely the exact computation of $I(a)$. This question, recently raised in \cite{bcs} in connection with the Hadamard conjecture, obviously stands above the Weingarten technology, and of the state-of-art of group integral computation in general. However, we do believe that the problem is actually quite tractable, and a ``magic formula'' for $I(a)$, in the case where $a$ is zero outside its upper left $2\times 2$ corner, is in preparation \cite{bc+}. There is of course an enormous work waiting to be done, for passing from the $2\times 2$ case to the $n\times n$ case, and it is our hope that the present results about $I(a)$ will be of help in this process.

The paper is organized as follows: 1-2 are preliminary sections, in 3-5 we state and prove the various Weingarten function results about $I(a)$, and 6-7 contain some results on the possible poles of $I(a)$, along with some concluding remarks.

\subsection*{Acknowledgements}

The present paper is part of a large computational project for polynomial group integrals, and it is a pleasure to thank B. Collins, S. Curran, E. Maurel-Segala, J. Novak, J.-M. Schlenker and P. Zinn-Justin for various discussions on the subject. The author was supported by the ANR grants ``Galoisint'' and ``Granma''.

\section{General considerations}

We are interested in the computation of arbitary polynomial integrals over the orthogonal group $O_n$. These are best introduced in a ``rectangular way'', as follows.

\begin{definition}
Associated to any matrix $a\in M_{p\times q}(\mathbb N)$ is the integral
$$I(a)=\int_{O_n}\prod_{i=1}^p\prod_{j=1}^qu_{ij}^{a_{ij}}\,du$$
with respect to the Haar measure of $O_n$, where $n\geq p,q$.
\end{definition}

We can of course complete our matrix with 0 values, as to always deal with square matrices, $a\in M_n(\mathbb N)$. However, the parameters $p,q$ are very useful, because they measure the ``complexity'' of the problem, as shown for instance by the result below.

Let $x!!=(x-1)(x-3)(x-5)\ldots$, with the product ending at $1$ or $2$.

\begin{theorem}
At $p=1$ we have the formula
$$I\begin{pmatrix}a_1&\ldots&a_q\end{pmatrix}=\varepsilon\cdot\frac{(n-1)!!a_1!!\ldots a_q!!}{(n+\Sigma a_i-1)!!}$$
where $\varepsilon=1$ if all $a_i$ are even, and $\varepsilon=0$ if not.
\end{theorem}

\begin{proof}
This follows from the fact that the first slice of $O_n$ is isomorphic to the real sphere $S^{n-1}$. Indeed, this gives the following formula:
$$I\begin{pmatrix}a_1&\ldots&a_q\end{pmatrix}=\int_{S^{n-1}}x_1^{a_1}\ldots x_q^{a_q}\,dx$$ 

This latter integral can be computed by using polar coordinates, and we obtain the formula in the statement. See e.g. \cite{bcs}. 
\end{proof}

Another instructive computation, as well of trigonometric nature, is the one at $n=2$. We have here the following result, which completely solves the problem in this case.

\begin{theorem}
At $n=2$ we have the formula
$$I\begin{pmatrix}a&b\\ c&d\end{pmatrix}=\varepsilon\cdot\frac{(a+d)!!(b+c)!!}{(a+b+c+d+1)!!}$$
where $\varepsilon=1$ if $a,b,c,d$ are even, $\varepsilon=-1$ is $a,b,c,d$ are odd, and $\varepsilon=0$ if not.
\end{theorem}

\begin{proof}
When computing the integral over $O_2$, we can restrict the integration to $SO_2=S^1$, then further restrict the integration to the first quadrant. We get: 
$$I\begin{pmatrix}a&b\\ c&d\end{pmatrix}
=\varepsilon\cdot\frac{2}{\pi}\int_0^{\pi/2}(\cos t)^{a+d}(\sin t)^{b+c}\,dt$$

This gives the formula in the statement.
\end{proof}

The above computations might tend to suggest that $I(a)$ always decomposes as a product of factorials. This is far from being true, but in the $2\times 2$ case it will be shown in \cite{bc+} that $I(a)$ decomposes as a (quite reasonable) sum of products of factorials.

\section{The Weingarten formula}

In this reminder of this paper we discuss the representation theory approach to the computation of $I(a)$. The technology available here is quite advanced, and will lead to a wealth of concrete results about $I(a)$, regarding its vanishing, theoretical poles, and asymptotic behaviour, and notably with a concrete formula for its first order term.

Our main tool is a combinatorial formula for the polynomial integrals over $O_n$, whose origins go back to Weingarten's paper \cite{wei}. The treatment given here follows the papers \cite{col}, \cite{csn}, \cite{cma}, where this formula was fully formalized, and its true power revealed.

The first remark is that the integrals over $O_n$ of arbitrary polynomial quantities of type $u_{i_1j_1}\ldots u_{i_sj_s}$ vanish, unless $s$ is even. In what follows, we focus on the $s=2k$ case.

Let us recall that the pseudo-inverse of a real symmetric matrix $G$ is the unique matrix $W$ satisfying $WGW=W$ and $GWG=G$. If $G$ is invertible, then $W=G^{-1}$.

\begin{theorem}
We have the Weingarten formula
$$\int_{O_n}u_{i_1j_1}\ldots u_{i_{2k}j_{2k}}\,du=\sum_{\pi,\sigma}\delta_\pi(i)\delta_\sigma(j)W_{kn}(\pi,\sigma)$$
where the various objects on the right are as follows:
\begin{enumerate}
\item The sum is over all pairings $\pi,\sigma$ of the set $\{1,\ldots,2k\}$.

\item The delta symbols are $0$ or $1$, depending on whether the indices $i=(i_1,\ldots,i_{2k})$ and $j=(j_1,\ldots,j_{2k})$ fit or not inside the corresponding pairings.

\item $W_{kn}$ is the pseudo-inverse of the matrix $G_{kn}(\pi,\sigma)=n^{loops(\pi,\sigma)}$, where we denote by $loops(\pi,\sigma)$ the number of loops obtained when superposing $\pi,\sigma$.
\end{enumerate}
\end{theorem}

\begin{proof}
Consider the matrix $P\in M_{n^{2k}}(\mathbb R)$ formed by all the integrals on the left, with $k$ fixed and $i,j$ varying. It follows from the general theory that $P$ is the orthogonal projection onto the space $Fix(u^{\otimes 2k})$ of fixed vectors of the $2k$-th tensor power of $u$.

By a well-known result of Brauer \cite{bra} the space $Fix(u^{\otimes 2k})$ is spanned by the vectors $\xi_\pi=\Sigma_i\delta_\pi(i)e_{i_1}\otimes\ldots\otimes e_{i_{2k}}$, with $\pi$ ranging over all the pairings of $\{1,\ldots,2k\}$. Now since the Gram matrix of these vectors is $<\xi_\pi,\xi_\sigma>=G_{kn}(\pi,\sigma)$, when computing $P$ we have to invert this Gram matrix, and we obtain the formula in the statement.
\end{proof}

We refer to \cite{col}, \cite{csn}, \cite{cma} for full details regarding the above proof. Let us also mention that the Gram matrix $G_{kn}$ is actually invertible for $n\geq k$. See \cite{csn}.

Our first task is to convert the above formula, by using our compact notation $I(a)$ for the polynomial integrals over $O_n$. We restrict attention to the case where $\Sigma a_{ij}$ is even.

\begin{theorem}
We have the Weingarten formula
$$I(a)=\sum_{\pi,\sigma}\delta_\pi(a_l)\delta_\sigma(a_r)W_{kn}(\pi,\sigma)$$
where $k=\Sigma a_{ij}/2$, and where the multi-indices $a_l/a_r$ are defined as follows:
\begin{enumerate}
\item Start with $a\in M_{p\times q}(\mathbb N)$, and replace each $ij$-entry by $a_{ij}$ copies of $i/j$.

\item Read this matrix in the usual way, as to get the multi-indices $a_l/a_r$.
\end{enumerate}
\end{theorem}

\begin{proof}
This is simply a reformulation of Theorem 2.1. Indeed, according to our definitions, the integral in the statement is given by:
$$I(a)=\int_{O_n}\underbrace{u_{11}\ldots u_{11}}_{a_{11}}\,\underbrace{u_{12}\ldots u_{12}}_{a_{12}}\,\ldots\,\underbrace{u_{pq}\ldots u_{pq}}_{a_{pq}}\,du$$

Thus what we have here is an integral as in Theorem 2.1, the multi-indices being:
\begin{eqnarray*}
a_l&=&(\underbrace{1\ldots 1}_{a_{11}}\,\underbrace{1\ldots 1}_{a_{12}}\,\ldots\,\underbrace{p\ldots p}_{a_{pq}})\\
a_r&=&(\underbrace{1\ldots 1}_{a_{11}}\,\underbrace{2\ldots 2}_{a_{12}}\,\ldots\,\underbrace{q\ldots q}_{a_{pq}})
\end{eqnarray*}

The result follows now from the Weingarten formula. 
\end{proof}

We are now in position of deriving a first general corollary from our study. This extends the vanishing results appearing in Theorem 1.2 and Theorem 1.3.

\begin{corollary}
We have $I(a)=0$, unless the matrix $a$ is ``admissible'', in the sense that all $p+q$ sums on its rows and columns are even numbers.
\end{corollary}

\begin{proof}
Observe first that the left multi-index associated to $a$ consists of $k_1=\Sigma a_{1j}$ copies of $1$, $k_2=\Sigma a_{2j}$ copies of $2$, and so on, up to $k_p=\Sigma a_{pj}$ copies of $p$. In the case where one of these numbers is odd we have $\delta_\pi(a)=0$ for any $\pi$, and this gives $I(a)=0$.

A similar argument with the right multi-index associated to $a$ shows that the sums on the columns of $a$ must be even as well, and we are done.
\end{proof}

\section{Basic asymptotic study}

A natural question now is whether the converse of Corollary 2.3 holds, and if so, the question of computing the sign of $I(a)$ appears as well. These are both quite subtle questions, and we begin our investigations with a $n\to\infty$ study.

The basic result here, known since \cite{col}, \cite{csn}, states that the Weingarten matrix is asymptotically diagonal, in the sense that we have $W_{kn}(\pi,\sigma)=n^{-k}(\delta_{\pi\sigma}+O(n^{-1/2}))$.

We present below a complete proof for this fact, by following a slighly improved method, which gives a better estimate, along with some concrete bounds.

\begin{theorem}
The Weingarten matrix is asymptotically diagonal, in the sense that:
$$W_{kn}(\pi,\sigma)=n^{-k}(\delta_{\pi\sigma}+O(n^{-1}))$$

\noindent Moreover, the $O(n^{-1})$ reminder is asymptotically smaller that $(2k/e)^kn^{-1}$.
\end{theorem}

\begin{proof}
It is convenient, for the purposes of this proof, to drop the indices $k,n$. We know that the Gram matrix is given by $G(\pi,\sigma)=n^{loops(\pi,\sigma)}$, so we have:
$$G(\pi,\sigma)=
\begin{cases}
n^k&{\rm for\ }\pi=\sigma\\
n,n^2,\ldots,n^{k-1}&{\rm for\ }\pi\neq\sigma
\end{cases}$$

Thus the Gram matrix is of the form $G=n^k(I+H)$, with $||H||_\infty\leq n^{-1}$. Now recall that for any $K\times K$ matrix $X$, we have the following lineup of standard inequalities:
$$||X||_\infty\leq||X||\leq||X||_2\leq K||X||_\infty$$

In the case of our matrix $H$, the size is $K=(2k)!!$, so we have $||H||\leq Kn^{-1}$. Now from $(I+H)^{-1}=I-H+H^2-H^3+\ldots$ we get $||I-(I+H)^{-1}||\leq ||H||/(1-||H||)$, so:
\begin{eqnarray*}
||I-n^kW||_\infty
&=&||I-(1+H)^{-1}||_\infty\\
&\leq&||I-(1+H)^{-1}||\\
&\leq&||H||/(1-||H||)\\
&\leq&Kn^{-1}/(1-Kn^{-1})\\
&=&K/(n-K)
\end{eqnarray*}

Together with the standard estimate $K\approx(2k/e)^k$, this gives the result.
\end{proof}

\begin{theorem}
We have the estimate
$$I(a)=n^{-k}\left(\prod_{i=1}^p\prod_{j=1}^qa_{ij}!!+O(n^{-1})\right)$$
when all $a_{ij}$ are even, and $I(a)=O(n^{-k-1})$ if not.
\end{theorem}

\begin{proof}
By using Theorem 2.2 and Theorem 3.1, we obtain the following estimate:
\begin{eqnarray*}
I(a)
&=&\sum_{\pi,\sigma}\delta_\pi(a_l)\delta_\sigma(a_r)W_{kn}(\pi,\sigma)\\
&=&n^{-k}\sum_{\pi,\sigma}\delta_\pi(a_l)\delta_\sigma(a_r)(\delta_{\pi\sigma}+O(n^{-1}))\\
&=&n^{-k}\left(\#\{\pi|\delta_\pi(a_l)=\delta_\pi(a_r)=1\}+O(n^{-1})\right)
\end{eqnarray*}

In order to count the partitions appearing in the set on the right, let us go back to the multi-indices $a_l,a_r$ described in the proof of Theorem 2.2. It is convenient to view both these multi-indices in a rectangular way, as follows: 
$$a_l=\begin{pmatrix}
\underbrace{1\ldots 1}_{a_{11}}&\ldots&\underbrace{1\ldots 1}_{a_{1q}}\\
\dots&\dots&\dots\\
\underbrace{p\ldots p}_{a_{p1}}&\ldots&\underbrace{p\ldots p}_{a_{pq}}$$
\end{pmatrix}
\quad\quad\quad
a_r=\begin{pmatrix}
\underbrace{1\ldots 1}_{a_{11}}&\ldots&\underbrace{q\ldots q}_{a_{1q}}\\
\dots&\dots&\dots\\
\underbrace{1\ldots 1}_{a_{p1}}&\ldots&\underbrace{p\ldots p}_{a_{pq}}$$
\end{pmatrix}$$

In other words, the multi-indices $a_l/a_r$ are now simply obtained from the matrix $a$ by ``dropping'' from each entry $a_{ij}$ a sequence of $a_{ij}$ numbers, all equal to $i/j$.

These two multi-indices, now in matrix form, have total length $2k=\Sigma a_{ij}$. We agree to view as well any pairing of $\{1,\ldots,2k\}$ in matrix form, by following the same convention.

With this picture, the pairings $\pi$ which contribute are simply those interconnecting sequences of indices ``dropped'' from the same $a_{ij}$, and this gives the following results:

(1) In the case where one of the entries $a_{ij}$ is odd, there is no pairing that can contribute to the leading term under consideration, so we have $I(a)=O(n^{-k-1})$, and we are done.

(2) In the case where all the entries $a_{ij}$ are even, the pairings that contribute to the leading term are those connecting points inside the $pq$ ``dropped'' sets, i.e. are made out of a pairing of $a_{11}$ points, a pairing of $a_{12}$ points, and so on, up to a pairing of $a_{pq}$ points. Now since an $x$-point set has $x!!$ pairings, this gives the formula in the statement. 
\end{proof}

\section{Geodesic expansion method}

We have seen how to compute the asymptotic behavior of $I(a)$, by using basic estimates on the Weingarten matrix. However, our results so far are effective only in the ``non-degenerate'' case, when all the entries of $a$ are even numbers. In the general case, the computation of the asymptotic sign of $I(a)$ requires a finer knowledge of the Weingarten matrix, notably with the exact computation of its leading term in $n^{-1}$.

In addition, the problem of computing the higher order terms, which are usually required for delicate applications of the Weingarten formula, appears as well.

We investigate these questions in this section and in the next one. We use a method of Collins \cite{col}, further processed by Collins-\'Sniady \cite{csn}. Let us begin with a key definition.

\begin{definition}
The Brauer space $D_k$ is defined as follows:
\begin{enumerate}
\item The points are the Brauer diagrams, i.e. the pairings of $\{1,2,\ldots,2k\}$. 

\item The distance function is given by $d(\pi,\sigma)=k-loops(\pi,\sigma)$.
\end{enumerate}
\end{definition}

It is indeed well-known, and elementary to check, that $d$ satisfies the usual axioms for a distance function. This actually comes from some general categorical properties of the Brauer diagrams, which are valid in much more general situations. See \cite{bco}, \cite{cur}.

The Brauer space, which will play an important role in what follows, is by definition a metric space having $(2k)!!=1.3.5\ldots(2k-1)$ points. An interesting question is to find a ``geometric'' realization of this space. This will be discussed later on.

The series expansion of the Weingarten function in terms of paths on the Brauer space was originally found by Collins in \cite{col} in the unitary case, then by Collins and \'Sniady \cite{csn} in the orthogonal case. We present below a slightly modified statement, along with a complete proof, by using a somewhat lighter formalism. 

\begin{theorem}
The Weingarten function $W_{kn}$ has a series expansion in $n^{-1}$ of the form
$$W_{kn}(\pi,\sigma)=n^{-k-d(\pi,\sigma)}\sum_{g=0}^\infty K_g(\pi,\sigma)n^{-g}$$
where the objects on the right are defined as follows:
\begin{enumerate}
\item A path from $\pi$ to $\sigma$ is a sequence $p=[\pi=\tau_0\neq\tau_1\neq\ldots\neq\tau_r=\sigma]$.

\item The signature of such a path is $+$ when $r$ is even, and $-$ when $r$ is odd.

\item The geodesicity defect of such a path is $g(p)=\Sigma_{i=1}^rd(\tau_{i-1},\tau_i)-d(\pi,\sigma)$.

\item $K_g$ counts the signed paths from $\pi$ to $\sigma$, with geodesicity defect $g$.
\end{enumerate} 
\end{theorem}

\begin{proof}
Let us go back to the proof of Theorem 3.1. We can write $G_{kn}=n^{-k}(I+H)$, and in terms of the Brauer space distance, the formula of $H$ is simply:
$$H(\pi,\sigma)=
\begin{cases}
0&{\rm for\ }\pi=\sigma\\
n^{-d(\pi,\sigma)}&{\rm for\ }\pi\neq\sigma
\end{cases}$$

Consider now the set $P_r(\pi,\sigma)$ of $r$-paths between $\pi$ and $\sigma$. According to the usual rule of matrix multiplication, the powers of $H$ are given by:
\begin{eqnarray*}
H^r(\pi,\sigma)
&=&\sum_{p\in P_r(\pi,\sigma)}H(\tau_0,\tau_1)\ldots H(\tau_{r-1},\tau_r)\\
&=&\sum_{p\in P_r(\pi,\sigma)}n^{-d(\pi,\sigma)-g(p)}
\end{eqnarray*}

Thus by using the formula $(1+H)^{-1}=1-H+H^2-H^3+\ldots$, we obtain:
\begin{eqnarray*}
W_{kn}(\pi,\sigma)
&=&n^{-k}\sum_{r=0}^\infty(-1)^rH^r(\pi,\sigma)\\
&=&n^{-k-d(\pi,\sigma)}\sum_{r=0}^\infty\sum_{p\in P_r(\pi,\sigma)}(-1)^rn^{-g(p)}
\end{eqnarray*}

Now by rearranging the various terms of the double sum according to their geodesicity defect $g=g(p)$, this gives the formula in the statement.
\end{proof}

For the $I(a)$ translation of the above result, it is convenient to use the total length of a path, defined as $d(p)=\Sigma_{i=1}^rd(\tau_{i-1},\tau_i)$. Observe that we have $d(p)=d(\pi,\sigma)+g(p)$.

\begin{theorem}
The integral $I(a)$ has a series expansion in $n^{-1}$ of the form
$$I(a)=n^{-k}\sum_{d=0}^\infty H_d(a)n^{-d}$$
where the coefficient on the right can be interpreted as follows:
\begin{enumerate}
\item Starting from $a\in M_{p\times q}(\mathbb N)$, construct the multi-indices $a_l,a_r$ as usual.

\item Call a path ``$a$-admissible'' if its endpoints satisfy $\delta_\pi(a_l)=1$ and $\delta_\sigma(a_r)=1$.

\item Then $H_d(a)$ counts all $a$-admissible signed paths in $D_k$, of total length $d$.
\end{enumerate}
\end{theorem}

\begin{proof}
We combine first Theorem 2.2 with Theorem 4.2:
\begin{eqnarray*}
I(a)
&=&\sum_{\pi,\sigma}\delta_\pi(a_l)\delta_\sigma(a_r)W_{kn}(\pi,\sigma)\\
&=&n^{-k}\sum_{\pi,\sigma}\delta_\pi(a_l)\delta_\sigma(a_r)\sum_{g=0}^\infty K_g(\pi,\sigma)n^{-d(\pi,\sigma)-g}
\end{eqnarray*}

Let us denote by $H_d(\pi,\sigma)$ the number of signed paths between $\pi$ and $\sigma$, of total length $d$. In terms of the new variable $d=d(\pi,\sigma)+g$, the above expression becomes:
\begin{eqnarray*}
I(a)
&=&n^{-k}\sum_{\pi,\sigma}\delta_\pi(a_l)\delta_\sigma(a_r)\sum_{d=0}^\infty H_d(\pi,\sigma)n^{-d}\\
&=&n^{-k}\sum_{d=0}^\infty\left(\sum_{\pi,\sigma}\delta_\pi(a_l)\delta_\sigma(a_r)H_d(\pi,\sigma)\right)n^{-d}
\end{eqnarray*}

We recognize in the middle the quantity $H_d(a)$, and this gives the result.
\end{proof}

\section{M\"obius function estimates}

In this section we derive some concrete consequences from the abstract results in the previous section. First, let us recall the following result, due to Collins and \'Sniady \cite{csn}.

\begin{theorem}
We have the estimate
$$W_{kn}(\pi,\sigma)=n^{-k-d(\pi,\sigma)}(\mu(\pi,\sigma)+O(n^{-1}))$$
where $\mu$ is the M\"obius function.
\end{theorem}

\begin{proof}
We know from Theorem 4.2 that we have the following estimate:
$$W_{kn}(\pi,\sigma)=n^{-k-d(\pi,\sigma)}(K_0(\pi,\sigma)+O(n^{-1}))$$

Now since one of the possible definitions of the M\"obius function is that this counts the signed geodesic paths, we have $K_0=\mu$, and we are done.
\end{proof}

It is probably interesting to note here that in some more general frameworks, e.g. when $O_n$ is replaced by its free version $O_n^+$, the straightforward extension of Theorem 5.1 doesn't seem to hold in full generality. In fact, the computation of the first order term of the generalized Weingarten functions is an open, interesting question. See \cite{bco}, \cite{cur}, \cite{bcu}.

Let us go back now to our integrals $I(a)$. The analogue of the above result of Collins and \'Sniady, fully generalizing Theorem 3.2, is as follows:

\begin{theorem}
We have the estimate
$$I(a)=n^{-k-e(a)}(\mu(a)+O(n^{-1}))$$
where the objects on the right are as follows:
\begin{enumerate}
\item $e(a)=\min\{d(\pi,\sigma)|\pi,\sigma\in D_k,\delta_\pi(a_l)=\delta_\sigma(a_r)=1\}$. 

\item $\mu(a)$ counts all $a$-admissible signed paths in $D_k$, of total length $e(a)$.
\end{enumerate} 
\end{theorem}

\begin{proof}
We know from Theorem 4.3 that we have an estimate of the following type:
$$I(a)=n^{-k-e}(H_e(a)+O(n^{-1}))$$

Here, according to the various notations in the previous section, $e\in\mathbb N$ is the smallest total length of an $a$-admissible path, and $H_e(a)$ counts all signed $a$-admissible paths of total length $e$. Now since the smallest total length of such a path is of course attained when the path is just a segment, we have $e=e(a)$ and $H_e(a)=\mu(a)$, and we are done.
\end{proof}

Summarizing, we have now a full description of the asymptotic behavior of $I(a)$. To illustrate how Theorem 5.2 applies, let us recover the results in Theorem 3.2:
\begin{enumerate}
\item Assume first that all entries of $a$ are even. In this case there is at least one partition $\pi$ such that $\delta_\pi(a_l)=\delta_\pi(a_r)=1$, so we have $e(a)=0$, and we recover the leading term $n^{-k}$ from Theorem 3.2. As for the coefficient $\mu(a)$, this counts the partitions $\pi$ such that $\delta_\pi(a_l)=\delta_\pi(a_r)=1$, so we fully recover Theorem 3.2. 

\item Assume now that $a$ has at least one odd entry. Here there is no $\pi$ as above, so we have $e(a)\geq 1$, and we recover the quantity $O(n^{-k-1})$ from Theorem 3.2.
\end{enumerate}

Regarding now the higher order terms, the situation is quite complicated. In fact, Theorem 4.3 is quite difficult to use as stated, as it was the case in fact with Theorem 4.2 as well. Let us point out, however, that a positive answer to the following question would definitely change the situation, by making these two results fully efficient:

\begin{question}
Is there a ``geometric'' model for the Brauer space?
\end{question}

In other words, what we are asking for is a better model, of geometric or perhaps analytic nature, which would convert the quite heavy counting of geodesics on the Brauer space into a kind of calculus computation. As far as we know, there is no idea yet here.

\section{Zonal function approach}

We have now a quite satisfactory picture of $\lim_{n\to\infty}I(a)$. However, the computation of $I(a)$ for fixed $n$ remains a quite subtle problem, as shown for instance by Theorem 1.3.

As explained in the previous section, the exact formula coming from the geodesic expansion is quite difficult to handle, in the lack of a better understanding of the Brauer space $D_k$. However, there have been many attempts for working out the combinatorics of the geodesic expansion. Let us mention here the pioneering work of Collins in the unitary case \cite{col}, further processed by Collins and \'Sniady in the orthogonal case \cite{csn}.

A remarkable advance on this subject came quite recently, in the preprint of Collins and Matsumoto \cite{cma}. In this section we briefly explain their new formula, and we work out the corresponding consequences regarding the integrals $I(a)$.

The formula of Collins and Matsumoto \cite{cma} (see also Zinn-Justin \cite{zin}) is as follows.

\begin{theorem}
We have the formula
$$W_{kn}(\pi,\sigma)=\frac{\sum_{\lambda\vdash k,\,l(\lambda)\leq k}\chi^{2\lambda}(1_k)w^\lambda(\pi^{-1}\sigma)}{(2k)!!\prod_{(i,j)\in\lambda}(n+2j-i-1)}$$
where the various objects on the right are as follows:
\begin{enumerate}
\item The sum is over all partitions of $\{1,\ldots,2k\}$ of length $l(\lambda)\leq k$.

\item $w^\lambda$ is the corresponding zonal spherical function of $(S_{2k},H_k)$.

\item $\chi^{2\lambda}$ is the character of $S_{2k}$ associated to $2\lambda=(2\lambda_1,2\lambda_2,\ldots)$.

\item The product is over all squares of the Young diagram of $\lambda$.
\end{enumerate}
\end{theorem}

It is of course possible to deduce from this a new a formula for $I(a)$, just by putting together the formulae in Theorem 2.2 and Theorem 6.1. However, there are probably several non-trivial simplifications that might appear when doing the sum over $\pi,\sigma$, and we do not know how the final statement about $I(a)$ should look like.

Instead, let us just record the following consequence.

\begin{corollary}
The possible poles of $I(a)$ can be at the numbers
$$-(k-1),-(k-2),\ldots,2k-1,2k$$
where $k\in\mathbb N$ associated to the admissible matrix $a\in M_{p\times q}(\mathbb N)$ is given by $k=\Sigma a_{ij}/2$.
\end{corollary}

\begin{proof}
We know from Theorem 2.2 that the possible poles of $I(a)$ can only come from those of the Weingarten function. On the other hand, Theorem 6.1 tells us that these latter poles are located at the numbers of the form $-2j+i+1$, with $(i,j)$ ranging over all possible squares of all possible Young diagrams, and this gives the result.
\end{proof}

As explained in \cite{cma}, the main feature of the zonal function approach is its remarkable efficiency in respect with numeric computations. It is our hope that the above results would be of help in producing, perhaps via a computer program, more ``data'' for the exact computation of $I(a)$. See the comments in the concluding section below. 

\section{Concluding remarks}

We have seen in this paper that the representation theory approach to the computation of polynomial integrals over $O_n$, originating from Weingarten's paper \cite{wei}, formulated and developed by Collins and \'Sniady in \cite{col}, \cite{csn}, further processed in \cite{cma}, \cite{mno}, \cite{zin}, and helped with some new ideas from \cite{bco}, \cite{cur}, \cite{dif}, leads to a number of concrete results on the integrals of type $I(a)$, which can be regarded as being the main objects of study.

Together with some explicit new computations, to be presented in our forthcoming paper \cite{bc+}, and based this time on purely geometric methods, these results give hope for the existence of a general explicit correspondence of type $a\to I(a)$, with the output $I(a)\in\mathbb Q(n)$ depending somehow ``arithmetically'' on the input $a\in M_{p\times q}(\mathbb N)$. 

We do not know if it is so, but we intend to intensively pursue the exploration of this correspondence, which would be of course extremely useful for any kind of mathematical area where the Weingarten formula has already proved to be useful. Note in particular that the present results provide a much better hope for the program launched in \cite{bcs}.

\end{document}